\documentclass[10pt,twoside]{amsart}
\usepackage[margin=0.92in]{geometry}
\usepackage[utf8]{inputenc}
\usepackage{amsmath}
\usepackage{mathtools}
\usepackage{graphicx} 
\graphicspath{ {./images/} }
\usepackage{wrapfig}
\usepackage{amsmath}
\usepackage{amsthm}
\usepackage{amsfonts}
\usepackage{amssymb}
\usepackage[all]{xy}
\usepackage{alltt}
\usepackage{enumitem}
\usepackage{stmaryrd}
\usepackage{comment}
\usepackage{xspace}
\usepackage{xcolor}
\usepackage{comment}
\usepackage{float}

\theoremstyle{plain}
\newtheorem{prop}{Proposition}[section]
\newtheorem{conjecture}[prop]{Conjecture}
\newtheorem{lem}[prop]{Lemma}

\newtheorem{thm}[prop]{Theorem}
\theoremstyle{definition}
\newtheorem*{defn}{Definition}

\newtheorem{cor}[prop]{Corollary}
\newtheorem{remark}[prop]{Remark}
\newtheorem{example}[prop]{Example}

\usepackage[
pdfauthor={ESYZ},
pdftitle={},
pdfstartview=XYZ,
bookmarks=true,
colorlinks=true,
linkcolor=blue,
urlcolor=blue,
citecolor=blue,
bookmarks=false,
linktocpage=true,
hyperindex=true
]{hyperref}

\title{Trigonometric functions in the $p$-norm}
\author{S. Chebolu}
\address{Department of Mathematics \\
Illinois State University \\
Normal, IL 61790, USA}
\email{schebol@ilstu.edu, abhatfi@ilstu.edu, rmklet1@ilstu.edu, cnmoor2@ilstu.edu, eaward4@ilstu.edu}

 \author{A. Hatfield}

\author{R. Klette}

\author{ C. Moore}

\author{E. Warden}

\date{\today}

\begin{document}

\thanks{The first author is supported by Simons Foundation: Collaboration Grant for Mathematicians (516354). }


\begin{abstract}
  Trigonometry is the study of circular functions,  which are functions defined on the unit circle $x^2+y^2 =1$, where distances are measured using the Euclidean norm. When distances are measured using the $L_p$-norm, we get generalized trigonometric functions. These are parametrizations of the unit $p$-circle $|x|^p+|y|^p =1$. Investigating these new functions leads to interesting connections involving double angle formulas, norms induced by inner products, Stirling numbers,  Bell polynomials,  Lagrange inversion, gamma functions, and generalized $\pi$ values.
\end{abstract}
\maketitle

\section{Introduction}
It is a well-known fact that trigonometric functions are periodic:  if $f(x)$ is any trigonometric function, then $f(x + 2\pi) = f(x)$ for all values of $x$ in the domain of $f$.  Therefore, it is natural to define trigonometric functions on the unit circle, where all multiples of $2\pi$ are identified when we wrap the real line onto the circle. Because of this definition, trigonometric functions are also called circular functions. In this setting, the trigonometric functions $\sin t$ and $\cos t$ are just the unit circle's parametrization with respect to arc length.

Recall that the unit circle is the locus of all points in the plane $\mathbb{R}^2$  that are at a distance of one unit from the origin, where distances are measured using the standard Euclidean norm: $\lVert \vec{x} \rVert = (x_1^2+x_2^2)^{1/2}$.   What if we switch to an $L_p$-norm: $\lVert \vec{x} \rVert_p = (|x_1|^p+|x_2|^p)^{1/p}$, $(p \ge 1 )$? We then get a new family of curves defined by the equations $|x|^p+|y|^p =1$. These are called unit $p$-circles and are shown in the figure below. Because these curves are in between a square and a circle, they are also called squircles.

\begin{figure}[h]
    \centering
    \includegraphics[width=0.3\textwidth]{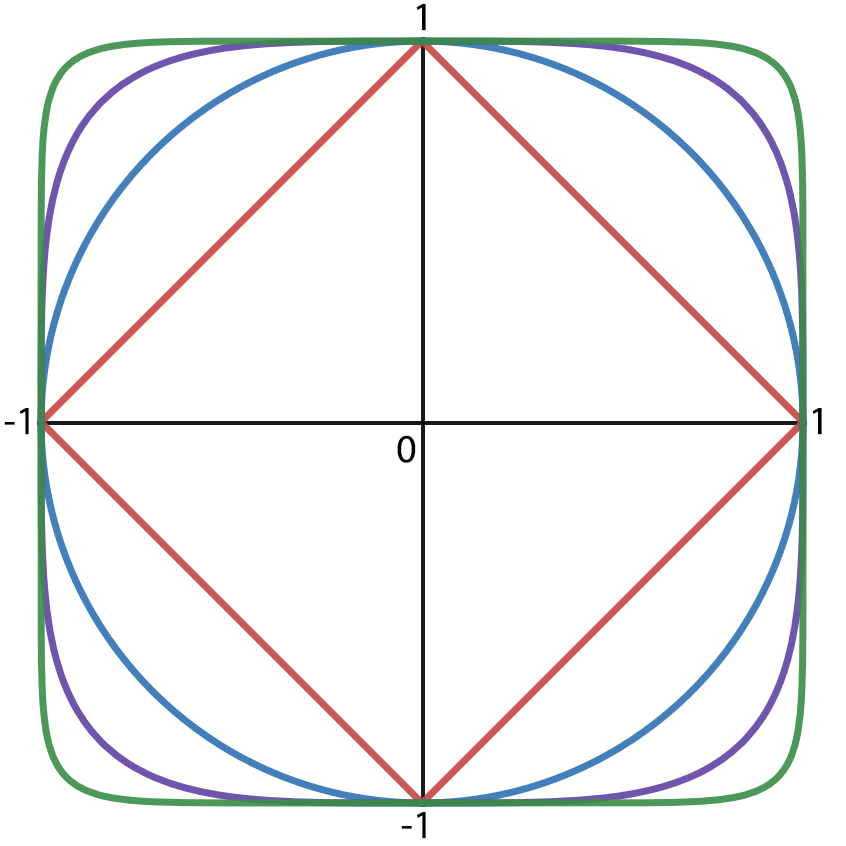}
    \caption{$p$-circles for $p = 1, 2, 4$, and $10$ from inside to outside, respectively }
\end{figure}

Can we parametrize these $p$-circles to get $p$-trigonometric functions $x = \sin_p t$ and $y= \cos_p t$ such that, when $p=2$, we recover the standard trigonometric functions? What properties and identities do these generalized trigonometric functions have? Can we do calculus over these curves? What can be said about the periods of these functions? How does the curvature change along a $p$-circle? What is the area it encloses? What are the rational points on $p$-circles? Note that for any $p \ge 1 $, $L_p$, as defined above, gives a norm, but this norm is induced by an inner product only when $p=2$ (\cite{SK}).
Therefore, $p=2$ is a special case of interest; however, all of the aforementioned questions are well defined for any $p \ge 1$.   The goal of this paper is to investigate these questions.  Our primary reference for this research is \cite{WP}. While we follow the general outline given in \cite{WP}, we also do some independent investigation. 

 There are at least three ways to generalize trigonometric functions. These correspond to 3 different \break parametrizations of the unit $p$-circle: areal, arc length, and angular.  It turns out that these three parametrizations are equivalent only when $p=2$! The parametrization we will be working with corresponds to the areal parametrization. Our investigation of these generalized trigonometric functions and their inverses led to several interesting connections involving double angle formulas, norms induced by inner products, Stirling numbers,  Bell polynomials,  Lagrange inversion, gamma functions,   and generalized $\pi$ values.

These $p$-trigonometric functions have several applications, specifically in design. Rather than using rounded rectangles, Apple uses $p$-circles for their icons, as the curvature continuity leads to a more sleek look, unifying the design of their hardware and icons \cite{AS}. Another design application can be found in squircular dinner plates, designed to allow a greater surface area for food while taking up the same amount of cabinet space as their circular counterparts \cite{PL}.

The paper is organized as follows. In Section \ref{ptrig}, we define $p$-trigonometric functions using a differential equations approach and derive some basic properties of these functions. We show that for any positive integer $k$, the well-known double angle formula for  $\sin(2x)$ holds for $\sin_k(2x)$ if and only if $k=2$.  In Section \ref{taylor}, we focus on the successive derivatives of $\sin_p(x)$. This revealed a connection between the coefficients of the terms in the derivatives and Stirling numbers of the first kind. We derive the Taylor series of $\sin_p^{-1}x$ using Newton's binomial series and then find the Taylor series of its inverse using Lagrange inversion theorem. It is shown that both $\sin_p x$ and $\sin_p^{-1} x$ are analytic functions at $x = 0$.  Our work gave rise to the concept of rigidity of functions, which deals with the simultaneous vanishing of the derivatives of a function and its inverse.  A generalization of $\pi$ for $p$-circles, $\pi_p$, and its properties are examined in Section \ref{pi_p} using  beta and gamma functions. Furthermore, we use a Monte Carlo method to compute $\pi_p$. In Section \ref{optimal}, we determine the value of $p$ for which the unit $p$-circle is halfway between the unit circle and the square that contains it from the lenses of area, perimeter, and curvature.
Rational points on $p$-circles are determined in Section \ref{rational}. We end the paper with some questions for future work in
Section \ref{conclusion}.

\vskip 3mm \noindent
\textbf{Acknowledgements:}  This paper is the outcome of the MAT 268 (Introduction to Undergraduate Research in Mathematics) course taught by the first author to the remaining authors at Illinois State University in Spring 2021. We want to thank the department of mathematics for providing us with the necessary resources for this research.  Discussions with Anindya Sen led to the notion of rigidity in Section \ref{taylor}.  Pisheng Ding raised several interesting questions and comments after reading this paper. We thank both of them for their interest and input in this paper. Finally, we are grateful to an anonymous referee for many comments and suggestions.

\section{$p$-trigonometric functions} \label{ptrig}
Unless stated otherwise, $p$ will denote a positive real number that is at least 1.

\subsection{Coupled Initial Value Problem}\label{CIVP} 
The standard trigonometric functions sine and cosine that parametrize the unit circle are famously coupled by the derivative relation 
$\sin't = \cos t,\, \cos't = -\sin t$. If we take $x(t) = \cos t$ and $y(t) = \sin t$, 
we see that the pair is one of many solutions to the system of differential equations
$$x'(t) = -y(t), \ \ \ y'(t) = x(t).$$ However, with the inclusion of the initial conditions 
$$x(0) = 1,\  y(0) = 0,$$ differential equation theory guarantees that the sine and cosine functions are, in fact, the only solutions
to this system \cite{BP}, better known as the Coupled Initial Value Problem (CIVP).

For  $p \geq 1$, a natural extension of the CIVP considers the functions $x(t),\, y(t)$ satisfying
$$x'(t) = -y(t)^{p-1}, \ y'(t) = x(t)^{p-1}, \ x(0) = 1, \ y(0) = 0.$$ 

The motivation for this extension comes from that fact that any functions $x(t)$ and $y(t)$ that satisfy the above CIVP parametrize the curve $x^p + y^p = 1$. This is seen by differentiating  $h(t) := x(t)^p + y(t)^p$ with respect to $t$, to get  $h'(t) = px(t)^{p-1}x'(t) + py(t)^{p-1}y'(t)$. Substituting $x'(t) = -y(t)^{p-1}, \ y'(t) = x(t)^{p-1}$, will show that $h'(t) = 0$. This means $h(t)$ is a constant function. Using the initial conditions, we can conclude that $h(t) = 1$, i.e., $x^p +y^p = 1$, as desired. 

Again, from the general theory of differential equations, the above CIVP has a unique solution. We can  define $\cos_pt = x(t)$
and $\sin_pt = y(t)$ as the unique solution to the generalized CIVP.  But these functions do not parametrize $p$-circles in general. For instance, when $p$ is an odd positive integer, these functions parametrize $p$-circles only in the first quadrant where $x$ and $y$ are both positive. To circumvent this issue, we restrict the domain of the solutions of the CIVP  and then extend them to functions on the real line using symmetry and periodicity. This is done in the next three subsections.

Once we have $\sin_p t$ and $\cos_p t$ in place, we may then define the other trigonometric functions $\tan_p t := \frac{\sin_p t}{\cos_p t},\, \csc_p t := \frac{1}{\sin_pt},\, \sec_pt := \frac{1}{\cos_pt}$, and $\cot_pt := \frac{1}{\tan_pt}$ such that the familiar inverse relations are maintained.

\subsection{Inverse $p$-trigonometric functions}\label{Inverse functions}
Starting with the equation $x=\sin_p{y}$, we use the CIVP to find $\sin_p^{-1}{x}$. Differentiating both sides with respect to $y$ and simplifying, we find:
\begin{eqnarray*}
\frac{dx}{dy} &=& \cos_p^{p-1}{y}\\
&=& (\cos_p^p{y})^{\frac{p-1}{p}}\\
&=& (1-\sin_p{y}^p)^{\frac{p-1}{p}} \\
&=& (1-x^p)^{\frac{p-1}{p}}.
\end{eqnarray*}
This is a separable differential equation. To solve it, we separate and integrate both sides. This gives:
\begin{eqnarray*}
\frac{dx}{dy} &=& (1-x^p)^{\frac{p-1}{p}}\\
\int{\frac{dx}{(1-x^p)^{\frac{p-1}{p}}}} &=& \int{dy} \\
\int_0^x{\frac{dt}{(1-t^p)^{\frac{p-1}{p}}}} &=& y = \sin_p^{-1}{x}.
\end{eqnarray*}
We can do the same for $x=\cos_p{y}$ to get
\begin{eqnarray*}
\cos_p^{-1}{x}= \int_x^1{\frac{dt}{(1-t^p)^{\frac{p-1}{p}}}}.
\end{eqnarray*}

\subsection{Areal parametrization of $p$-circles}
The unit circle has a useful property that a sector with angle measure $\theta$ in radians has an area of $\theta/2$. We can use this property to find sine and cosine in terms of area where $x=\cos(2a)$, $y=\sin(2a)$, and $a$ is the area of the sector made by the points $(1,0)$ and $(x,y)$. It is then natural to ask if this property extends to all $p$-circles.

\begin{figure}
    \centering
    \includegraphics[width=0.3\textwidth]{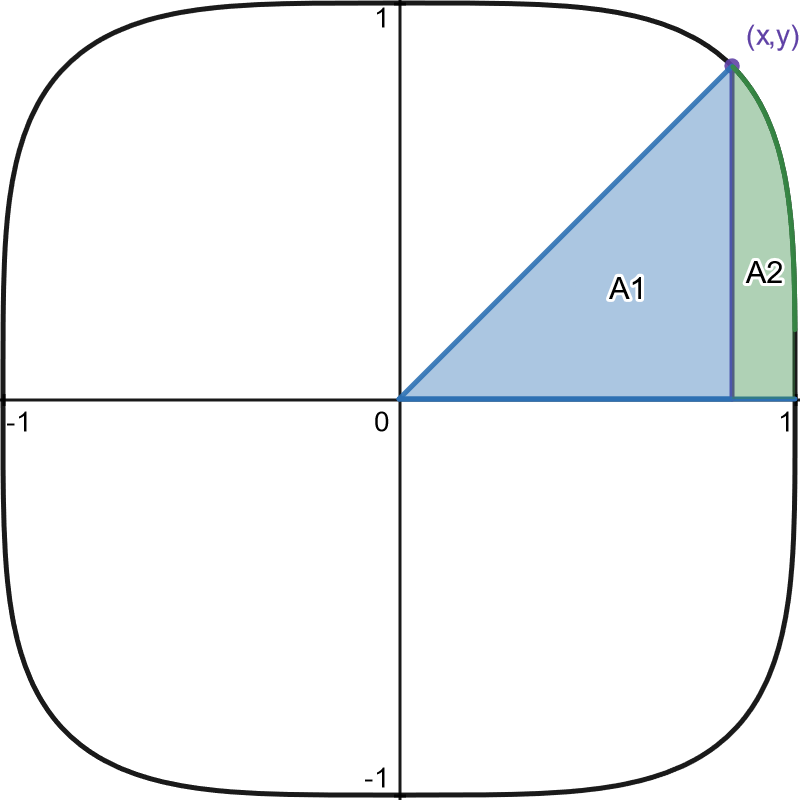}
    \caption{Area of a $p$-sector}
    \label{fig:areaOfSector}
\end{figure}

\begin{prop}\label{areapara}
Let  $(x,y)$ be a point in the first quadrant of the unit $p$-circle, and $a$ be the area of the sector made by the points $(1,0)$ and $(x,y)$. It holds that $x=\cos_p(2a)$ and $y=\sin_p(2a)$.
\end{prop}

\begin{proof}
This argument is in the spirit of Levin \cite{LA}. Working in the first quadrant, the area of the sector in a $p$-circle can be given by the area of $A1+A2$ as denoted in Figure \ref{fig:areaOfSector}. This can be given by $a=\frac{1}{2}x(1-x^p)^{\frac{1}{p}}+\int_x^1{(1-t^p)^{\frac{1}{p}}dt}$. We can differentiate both sides with respect to $x$ and simplify to get the following:
\begin{eqnarray*}
\frac{da}{dx} &=& \frac{1}{2}\left((1-x^p)^{\frac{1}{p}}+x\frac{1}{p}(1-x^p)^{\frac{1}{p}-1}(-px^{p-1})\right)-(1-x^p)^{\frac{1}{p}}\\
    &=&(1-x^p)^{\frac{1}{p}}\left(\frac{1}{2}-\frac{x^p}{2}(1-x^p)^{-1}-1 \right)\\
    &=&(1-x^p)^{\frac{1}{p}}\left(\frac{(1-x^p)-x^p-2(1-x^p)}{2(1-x^p)} \right)\\
    &=&-\frac{(1-x^p)^{\frac{1}{p}-1}}{2}.
\end{eqnarray*}

Using the fundamental theorem of calculus, we can write this as $a=\int_x^1{\frac{(1-t^p)^{\frac{1}{p}-1}}{2}dt}+c$. When $a=0$ and $x=1$, we get $c=0$. From here, we can conclude that $a=\frac{1}{2}\arccos_p{x}$. Solving for $x$ gives $x=\cos_p{(2a)}$.

We can do the same thing in terms of $y$ to get $a=\int_0^y{\frac{(1-t^p)^{\frac{1}{p}-1}}{2}dt}$+c. When $a=0$ and $y=0$, we get $c=0$. From here, this equation has been shown to be $a=\frac{1}{2}\arcsin_p{y}$ and thus $y=\sin_p{(2a)}$. As such, this shows that this property does extend to all unit $p$-circles.
\end{proof}

\subsection{Definition and Graphs of $\sin_px$ and $\cos_px$} \label{sec:extension}

To generalize the formula $\pi/2= \sin^{-1} (1)$, we first  set $\pi_p/2 :=  \sin_p^{-1} (1)$. 
  Since we have shown that the $p$-trigonometric functions can be parametrized by area, we can now extend then to functions defined on the entire real line as follows. We first restrict them to $[0, \sin_p^{-1} (1)] =[0 , \pi_p/2]$ and then  extend  the domain to $[0, 2\pi_p]$ using symmetry:
    \[
\sin_p t :=  \begin{cases}
    \sin_p(\pi_p - t)  & \pi_p/2 < t \le \pi_p, \\
    -\sin_p(2\pi_p - t)  & \pi_p < t < 2\pi_p.
    \end{cases}
    \]
We then periodically extend that it  $(-\infty,\infty)$ 
 by setting $\sin_p(t+2\pi_pk) = \sin_p(t)$ for any integer $k$. The definition of $\cos_p(t)$ is similar. The resulting graphs are shown below.

    \begin{figure}[H]
      \centering
      \includegraphics[width=1\textwidth]{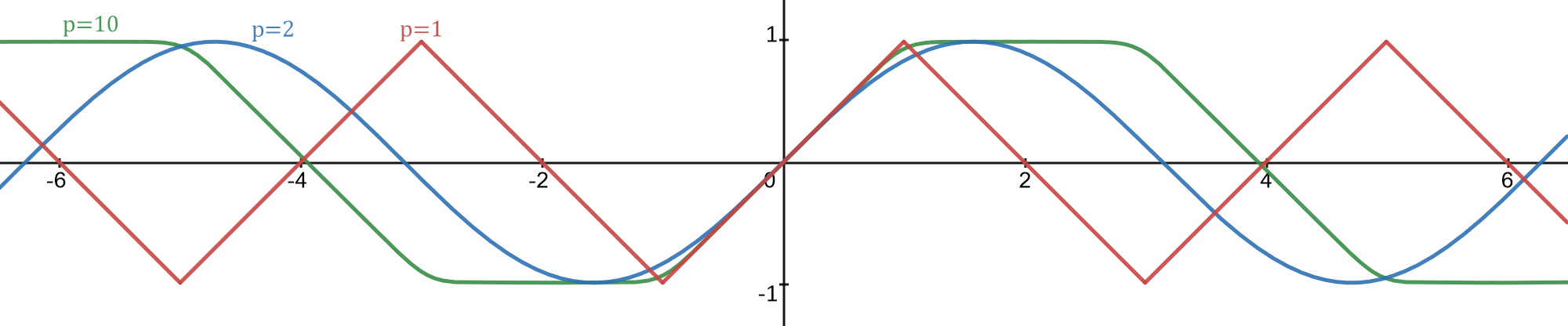}
      \caption{Graph of $\sin_p{x}$ for $p=1,2,$ and $10$}
      \label{fig:sinp}
    \end{figure}
    \begin{figure}[H]
      \centering
      \includegraphics[width=1\textwidth]{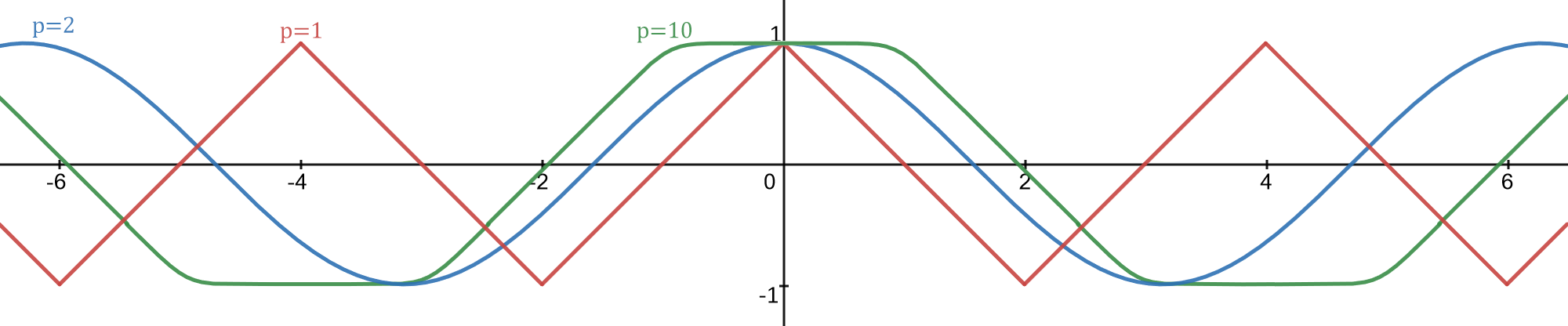}
      \caption{Graph of $\cos_p{x}$ for $p=1,2,$ and $10$}
      \label{fig:cosp}
    \end{figure}

\subsection{Trigonometric identities} 
While we may have defined the generalized CIVP in a manner similar to the original, there is no guarantee that $\sin_p t,\, \cos_p t$ thus defined satisfy familiar trigonometric properties and identities. In this section, we explore a few identities of the $p$-trigonometric functions.

\begin{lem}[$p$-Pythagorean Equation]\label{p-pythag} \cite[p. 268]{WP}
The functions $\sin_p t,\, \cos_p t$  satisfy $|\sin_pt|^p + |\cos_pt|^p = 1$ for all real $t$. 
\end{lem}

\begin{proof}
This is clear from the definition of these functions using the CIVP and extension using symmetry and periodicity.
\end{proof}

It is clear from the $p$-Pythagorean Equation that the functions $\sin_p t,\, \cos_pt$ are bounded and
$|\sin_pt| \leq 1,\, |\cos_pt| \leq 1$. Dividing all terms of the $p$-Pythagorean equation  by $|\sin_pt|^p$ and $|\cos_pt|^p$ gives the identities
$1 + |\cot_pt|^p = |\csc_pt|^p$ and $|\tan_pt|^p + 1 = |\sec_pt|^p$, respectively.

\begin{lem}
    The functions $\sin_pt$ and $\cos_pt$ are odd and even, respectively.
\end{lem}

\begin{proof}
The functions $\alpha(t) := -\sin_p(-t)$ and $\beta(t) = \cos_p(-t)$
satisfy $\alpha'(t) = -(-\sin_p'(-t)) = \cos_p^{p-1}(-t) = \beta(t)^{p-1}$ and 
$\beta'(t) = -\cos'_p(-t) = \sin_p^{p-1}(-t) = -\alpha(t)^{p-1}$. Note that 
$\alpha(0) = -\sin_p(0) = 0$ and $\beta(0) = \cos_p(0) = 1$; thus, the functions
$\alpha, \beta$ satisfy the generalized CIVP. Then, by the uniqueness of solutions, we must have
$\sin_pt = -\sin_p(-t)$ and $\cos_pt = \cos_p(-t)$. 
\end{proof}

However, not all standard $2$-trigonometric identities are satisfied. For instance, we show that for positive integer values of $p$,  $\sin_p(2t) = 2\sin_pt\cos_pt$ is satisfied if and only if $p = 2$.  
A double angle formula for generalized trigonometric functions is still sought after \cite{ED, SS}.

\begin{prop}\label{doubleangle}
 Let $k \in \mathbb{Z}_+$. Then $\sin_k(2t) = 2\sin_k(t)\cos_k(t)$ if and only if $k = 2$. 
\end{prop}

\begin{proof}
The desired identity is well known for $k = 2$. We suppose the identity holds for $k \ge 1$ and show that $k$ must be 2.
We consider the cases $k = 1$ and $k >  1$ separately. For $k = 1$, we note that the CIVP gives the unique solution
$\sin_1(t) = t$ and $\cos_1(t) = 1-t$. Then $\sin_1(2t) = 2t \ne 2t(1-t) = 2\sin_1 t \cos_1 t$. If
$k > 1$, then by Lemma \ref{p-pythag}, the functions $\sin_k$ and $\cos_k$ satisfy 
$|\sin_k t|^k + |\cos_k t|^k = 1$. By the Intermediate Value Theorem, there exists some $t_0$ in $[0, \pi_p/2]$
 such that $\sin_k(t_0) = \cos_k(t_0)$. As $t_0 \geq 0$, the substitution
$\sin_k(t_0) = \cos_k(t_0)$ into the $p$-Pythagorean identity gives 
$2\sin_k^k (t_0) = 1$, therefore $\sin_k^k (t_0) = \frac{1}{2}$. Then, by the assumption that $\sin_k(2t) = 2\sin_k t \cos_k t$ is satisfied for all $t$, we may
raise all terms to the power $k$ and evaluate at the point $t_0$ to obtain
$\sin_k^k(2t_0) = 2^k (\frac{1}{2}) \frac{1}{2} = 2^{k-2}$. Since $\sin_k^k t$ is bounded above by $1$, we obtain
$2^{k-2} \leq 1$, which implies that $k \leq 2$. Together with the assumption that $k > 1$, we obtain 
$k = 2$. 
\end{proof}

It is known that the $L_p$ norm is induced by an inner product if and only if $p = 2$ \cite{SK}. Then together with 
Proposition \ref{doubleangle}, we make the following remark. 

\begin{remark}
The following are equivalent for $k \in \mathbb{Z}_+$:
\begin{itemize}
    \item $L_k$ is a norm induced by an inner product,
    \item $\sin_k(2t) = 2\sin_k t \cos_k t$, and
    \item $k = 2$.
\end{itemize}
\end{remark}

\section{Taylor Series} \label{taylor}
Now that we have defined $p$-trigonometric functions and their derivatives by the CIVP, it is natural to study the higher derivatives of these functions.  We begin by observing that for any $p > 1$, all the successive derivatives of $\sin_p x$ and $\cos_p x$ are defined for all values of $x$. In this section, we provide an algorithm for differentiating these functions, demonstrate some patterns and connections present in their successive derivatives, and formulate the Taylor series for $\sin_p^{-1}x$ and $\sin_p x$. The Taylor series representations of these functions provide a tool to express all of the derivatives of the $p$-trigonometric functions in one formula.
\subsection{Higher derivatives and the bracket notation}
Because of the simplicity and utility of the closed formulas for differentiation of $\sin_2x,\, \cos_2x$, it is natural to wonder about higher derivatives of $\sin_px,\, \cos_px$. We find these higher derivatives by utilizing the definition given by the CIVP in Section \ref{CIVP}. However, these derivatives become complex rather quickly. To help address this, we introduce a notation that will be used throughout this section in relation to higher derivatives of these $p$-trigonometric functions: $[m,n]_p:=\cos_p^m(x)\sin_p^n(x)$.
\begin{lem}\label{bracketNotaiton}
    The derivative of $\cos_p^m(x)\sin_p^n(x)$ satisfies $ \frac{d}{dx}[m,n]_p=-m[m-1,n+p-1]_p+n[m+p-1,n-1]_p$.
\end{lem}

\begin{proof} Applying the standard rules of differentiation, we get the following. 
\begin{align*}
 \frac{d}{dx} [m,n]_p &= 
     \frac{d}{dx} (\cos_p^m(x)\sin_p^n(x)) \\
   & = -m\cos_p^{m-1}(x)\sin_p^{p-1}(x)\sin_p^n(x)+\cos_p^m(x)\cdot n\sin_p^{n-1}(x)\cos_p^{p-1}(x) \\
    &=-m\cos_p^{m-1}(x)\sin_p^{n+p-1}(x)+n\cos_p^{m+p-1}(x)\sin_p^{n-1}(x). \\
     &=-m[m-1,n+p-1]_p+n[m+p-1,n-1]_p. \qedhere
\end{align*}
\end{proof}

Although we do not have a closed formula for finding derivatives of these functions, Lemma \ref{bracketNotaiton} serves as a recursive algorithm for computing successive derivatives, as demonstrated in the following example.

\begin{example}\label{derivativeExample}
    Lemma \ref{bracketNotaiton} can be iteratively applied to $\sin_px$ to find the first few derivatives:
    \begin{align*}
        \sin_p x&=[0,1]_p \\
        \frac{d}{dx}  \sin_p x&=0+1[p-1,0]_p\\
        \frac{d^2}{dx^2} \sin_p x&=0+(-p+1)[p-2,p-1]_p+0 \\
        \frac{d^3}{dx^3}  \sin_px&=0+-(p-1)(-(p-2)[p-3,2p-2]_p+(p-1)[2p-3,p-2]_p)+0 \\
        &=0+(p^2-3p+2)[p-3,2p-2]_p+(-p^2+2p-1)[2p-3,p-2]_p+0.
    \end{align*}
\end{example}
There seems to be no clear pattern that arises from these derivatives like there is for $\sin x$. However, in the next subsection, we will see one pattern in the coefficients of the first terms of these derivatives.

\subsection{Connection to Stirling numbers}
For any variable $x$ and a non-negative integer $n$, the falling factorial is defined as follows.
\[(x)_n:= \begin{cases}
1 & \text{ if } n =0, \\
x(x-1)(x-2) \cdots (x-n+1) & \text{ if } n \ge 1 .
\end{cases}\]

For $n \ge 1$, $(x)_n$ is a non-constant polynomial of degree $n$ whose coefficients are the Stirling numbers of the first kind. More precisely, we set: 
\[ (x)_n  = \sum_{k=1}^n s(n, k) x^k.\]

We will now show a connection between the successive derivatives of $\sin_p x$ and Stirling numbers.
Building a tower from the coefficients in Example \ref{derivativeExample}, we get: 
\begin{center}
    1\\
    0 \ \ \verb"|" \ \ \underline{1}\\
    0 \ \ \verb"|" \ \ \underline{$-p+1$} \ \ \verb"|" \ \ 0\\
    0 \ \ \verb"|" \ \ \underline{$p^2-3p+2$} \ \ \verb"|" \ \ $-p^2+2p-1$ \ \ \verb"|" \ \ 0
\end{center}
We observed that  the coefficients of the polynomials in the second column (underlined) can be expressed using Stirling numbers of the first kind $s(n,k)$. For instance, corresponding to the polynomial $p^2-3p+2$ (corresponding to  the 3rd derivative of $\sin_p x$), we have $s(3, 3) = 1$, $s(3, 2) = -3$ and $s(3, 1) = 2$. To prove this, we need the following lemma.

\begin{lem}
For any $n \ge 1$, the first term of $\frac{d^n}{dx^n} (\sin_p(x))$ is given by 
\[ (-1)^{n-1} (p-1)_{n-1} [p-n, (n-1)(p-1)]_p.\]
\end{lem}

\begin{proof}
We prove this using mathematical induction. For $n=1$, $\frac{d}{dx} (\sin_p(x)) = \cos_p^{p-1} (x) = 1 [p-1, 0]_p$, which agrees with the answer obtained with $n=1$ in the given expression. Having proved the base case, let us assume that the result is true for $n=k$. Differentiating the first term of $\frac{d^k}{dx^k} (\sin_p(x))$ using the chain rule, and only picking the first term of the resulting expression will give us
\[ (-1)^{k-1} (p-1)_{k-1} \left( (-1) (p-k) [p-(k+1), k (p-1)] \right) = (-1)^k (p-1)_{k-1} (p-1) [p-(k+1), k(p-1)].\]
The recursive nature of the falling factorial tells us that $(p-1)_{k-1} (p-1) = (p-1)_k$. This shows that the first term of $\frac{d^{k+1}}{dx^{k+1}} (\sin_p(x))$ is given by $(-1)^k(p-1)_k  [p-(k+1), k(p-1)]$. By the principle of mathematical induction, the result is true for all $n \ge 1$.
\end{proof}

The connection to Stirling numbers and the successive derivatives of the $\sin_p(x)$ is now clear. Simplifying the   coefficient of the first term of $\frac{d^n}{dx^n} (\sin_p(x))$ obtained from the above lemma gives:
\[ (-1)^{n-1} (p-1)_{n-1} = (-1)^{n-1} \frac{(p)_{n}}{p} = \frac{(-1)^{n-1}}{p}  \sum_{k=1}^n s(n, k) p^k.\]

\subsection{Newton's binomial series} Let $p$ be any integer that is greater than 1.
As the previous section demonstrates, finding a formula for the successive derivatives of $\sin_px$ to compute its Taylor series is complicated. Instead, we examine $\sin_p^{-1}x$, whose Taylor series at $x=0$  is more manageable, and use this to find the Taylor series of $\sin_px$ at $x= 0$ through the Lagrange inversion theorem. To do this, we apply Newton's binomial series to derive the Taylor series of $\sin_p^{-1}x$. Newton's binomial series tells us the following for any exponent $a$ and $|x| < 1$:
\begin{align*}
    (1-x)^{-a}&=1+ax+\frac{a(a+1)}{2!}x^2+\frac{a(a+1)(a+2)}{3!}x^3+... \\
    &=\sum_{k=0}^\infty\frac{a^{(k)}x^k}{k!},
\end{align*}
where $a^{(k)}=a(a+1)(a+2)\cdots(a+k-1)$ is the rising factorial \cite[p. 742]{JS}. Note that, by convention, $a^{(0)} = 1$.
\begin{prop} \label{taylor-arcsin}
    We can express $\sin_p^{-1}x$ as the following Taylor series:
    \begin{equation*}
        \sin_p^{-1}x=\sum_{k=0}^\infty\left(\frac{p-1}{p}\right)^{(k)}\frac{x^{kp+1}}{k!(kp+1)}.
    \end{equation*}
\end{prop}

\begin{proof}
    Beginning with the integral form of $\sin_p^{-1}x$ derived in Section \ref{Inverse functions}, we apply Newton's binomial series:
\begin{align*}
    \sin_p^{-1}x&=\int_0^x(1-t^p)^{-(\frac{p-1}{p})}dt \\
    &=\int_0^x\left(\sum_{k=0}^\infty\left(\frac{p-1}{p}\right)^{(k)}\frac{t^{kp}}{k!} \right) dt.
\end{align*}
Power series have the property that they can be integrated term by term within the interval of convergence. Thus, when we integrate and apply the fundamental theorem of calculus, the result follows.
\end{proof}

\begin{example} \label{talorexample}
    Applying Proposition \ref{taylor-arcsin} for $p=2$ gives the  following well-known result:
    \begin{equation*}
        \sin_2^{-1}x=x+\frac{1}{6}x^3+\frac{3}{40}x^5+\frac{5}{112}x^7+ \cdots +{2n \choose n}\frac{x^{2n+1}}{2^{2n}(2n+1)}+ \cdots.
    \end{equation*}
    Similarly, when $p=4$, we get the first few terms as follows:
    \begin{equation*}
        \sin_4^{-1}x=x+\frac{3}{20}x^5+\frac{7}{96}x^9+\frac{77}{1664}x^{13}+ \cdots.
    \end{equation*}
\end{example}
It would be helpful to have a closed-form solution for these higher derivatives. In the next section, we introduce some tools and discuss what this will look like.

\subsection{$\sin_p^{-1} x$ through the gamma function}
We now introduce a special function to shed light on $\sin_p^{-1} x$.  The gamma function, $\Gamma(z)$, is defined as $\Gamma(z)=\int_{0}^{\infty} e^{-t}t^{z-1} \,dt$, for $z > 0$. This converges for any real number $z > 0$, and it 
is an extension of the factorial function: $\Gamma(n) = (n-1)!$. It is well-known that $\Gamma(1/2) = \sqrt{\pi}$.  Two important properties of the gamma function are 
\[\Gamma(x+1) = x \Gamma(x) \ \ \text{ and } \ \ \Gamma(x) \Gamma(1-x) = \frac{\pi}{ \sin(\pi x)}.\]
\\

Using the gamma function, for any integer $p > 1$, we can further simplify the Taylor series for $\sin^{-1}_p(x)$ as follows.  
We begin by the formula from Proposition \ref{taylor-arcsin} which states that 
 \begin{equation*}
        \sin_p^{-1}x=\sum_{k=0}^\infty\left(\frac{p-1}{p}\right)^{(k)}\frac{x^{kp+1}}{k!(kp+1)}.
    \end{equation*}
Then we have the following:    
\begin{align*}
        \sin_p^{-1}x & = \sum_{k=0}^\infty\left(\frac{p-1}{p}\right)^{(k)}\frac{x^{kp+1}}{k!(kp+1)}   \\
        & = \sum_{k=0}^\infty\left(1-\frac{1}{p}\right)^{(k)}\frac{x^{kp+1}}{k!(kp+1)}  \\
        & = \sum_{k=0}^\infty \left(1-\frac{1}{p}\right) \left(2-\frac{1}{p}\right)\left(3-\frac{1}{p}\right) \cdots \left(k-\frac{1}{p}\right)\frac{x^{kp+1}}{k!(kp+1)}  \\
        & = \sum_{k=0}^\infty\Gamma\left(1- \frac{1}{p}\right) \left(1-\frac{1}{p}\right) \left(2-\frac{1}{p}\right)\left(3-\frac{1}{p}\right) \cdots \left(k-\frac{1}{p}\right) \frac{1}{\Gamma(1- \frac{1}{p}) }\frac{x^{kp+1}}{k!(kp+1)} \\
        & = \sum_{k=0}^\infty\Gamma\left(2- \frac{1}{p}\right)  \left(2-\frac{1}{p}\right)\left(3-\frac{1}{p}\right) \cdots \left(k-\frac{1}{p}\right) \frac{1}{\Gamma(1- \frac{1}{p}) }\frac{x^{kp+1}}{k!(kp+1)} \\
  & \ \  \vdots \\
         & = \sum_{k=0}^\infty\Gamma\left(k- \frac{1}{p}\right)  \frac{1}{\Gamma(1- \frac{1}{p}) }\frac{x^{kp+1}}{k!(kp+1)} \\
         & = \sum_{k=0}^\infty \frac{\Gamma\left(k- \frac{1}{p}\right)  }{\Gamma\left(1 - \frac{1}{p}\right)}\frac{x^{kp+1}}{k!(kp+1)}. 
\end{align*}

\begin{thm} \label{thm1}
Let $n > 1$ be a positive integer. Then for any positive integer $l$,  let $k$ and $r$ be the integers given by the division algorithm: $l = nk + r$ where $k \ge 0$ and  $0 \le r \le  n-1$. Then  we have
\[ \left(\frac{d^l}{dx^l} \sin_n^{-1}(x) \,\right) \biggr \rvert_{x=0} = \begin{cases}
\frac{\Gamma\left(k- \frac{1}{n}\right)  }{\Gamma\left(1 - \frac{1}{n}\right)}\frac{(kn)!}{k!},  &  \text{ if } r =  1,  \\
 0,  &  \text{ if } r \ne  1. 
\end{cases}
\]
\end{thm}

\begin{proof}
The Taylor series for $\sin^{-1}_n(x)$ at $x = 0$ has the form 
\[\sin^{-1}_n(x) = h(x) =  \sum_{m=0}^{\infty} \frac{h^{(m)}(0)}{m!} x^m.  \]
On the other hand, from the above calculation, we know that 
\[\sin^{-1}_n(x) = \sum_{k=0}^\infty \frac{\Gamma\left(k- \frac{1}{n}\right)  }{\Gamma\left(1 - \frac{1}{n}\right)}\frac{x^{kn+1}}{k!(kn+1)} = 
 \sum_{k=0}^\infty \left(\frac{\Gamma\left(k- \frac{1}{n}\right)  }{\Gamma\left(1 - \frac{1}{n}\right)}\frac{(kn)!}{k!} \right) \frac{1}{(kn+1)!} x^{kn+1}.\]
Equating the coefficients of like-powers of $x$ in both these series, we get the theorem.
\end{proof}

Now that we have derived the Taylor series of $\sin_p^{-1}x$, we can apply the Lagrange inversion theorem as outlined in the next section.

\subsection{Lagrange inversion}\label{lagrange}  A function $z = f(w)$ is said to be analytic at $c$ if it is infinitely differentiable at $c$ and if the Taylor series for $f(w)$ at $w=c$ converges to $f(w)$ for all $w$ in a neighborhood of $c$.

For an equation $z=f(w)$, where $f$ is analytic at $c$ and $f'(c)\neq 0$, the Lagrange inversion theorem can be used to find the equation's inverse, $w=g(z)$, in a neighborhood of $0$. This inverse is given by the formula \cite[Chapter 3]{MA}:
\begin{align*}
    g(z)&=c+\sum_{n=1}^\infty g_n\frac{(z-f(c))^n}{n!},\ \text{ where} \\
    g_n&=\lim_{w\to c}\frac{d^{n-1}}{dw^{n-1}} \left[\left(\frac{w-c}{f(w)-f(c)}\right)^n \right].
\end{align*}
    For power series, this theorem takes a slightly different form. Specifically, when $f$ and $g$ are formal power series expressed as
    \begin{equation*}
        f(w)=\sum_{k=0}^\infty f_k\frac{w^k}{k!}\quad \text{and}\quad g(z)=\sum_{k=0}^\infty g_k\frac{z^k}{k!},
    \end{equation*}
    with $f_0=0$ and $f_1\neq 0$, applying the Lagrange inversion theorem gives us the following \cite{MA}:
    \begin{align*}
        &g(z)=c+\sum_{n=1}^\infty g_n\frac{(z-f(c))^n}{n!},\quad \text{with} \\
        &g_n=\frac{1}{f_1^n}\sum_{k=1}^{n-1}(-1)^kn^{(k)}B_{n-1,k}(\hat{f}_1,\hat{f}_2,...,\hat{f}_{n-k}),\quad n\geq2,\quad \text{where} \\
        &\hat{f}_k=\frac{f_{k+1}}{(k+1)f_1}, \quad g_1=\frac{1}{f_1}, \quad n^{(k)}=n(n+1)\cdots(n+k-1)\text{, and} \\
        B_{n,k}(x_1,&x_2,...,x_{n-k+1})= \sum_ {} \frac{n!}{j_1!j_2!...j_{n-k+1}!}\left (  \frac{x_1}{1!}   \right )^{j_1} \left (  \frac{x_2}{2!}   \right )^{j_2} ... \left (  \frac{x_{n-k+1}}{(n-k+1)!}   \right )^{j_{n-k+1}},
    \end{align*}
    where this sum is taken over all sequences $j_1,j_2,j_3,..., j_{n-k+1}$ of non-negative integers that satisfy $j_1+j_2+...+j_{n-k+1}=k$ and
$j_1+2j_2+3j_3+...+(n-k+1)j_{n-k+1}=n$. These are the Bell polynomials.

    The Taylor series expansion of $\sin_p x$ is obtained when the above theorem is applied to
    \begin{equation*}
        \sin_p^{-1}x =\sum_{k=0}^\infty \left(\frac{p-1}{p}\right)^{(k)} \frac{x^{kp+1}}{k!(kp+1)},
    \end{equation*}
    which was derived in the previous section. We are able to apply this theorem to $\sin_p^{-1}x$, as it meets the initial conditions given above: $f_0=0$ and $f_1\neq 0$.

    \begin{example}
        When $p=2$, we can apply Lagrange Inversion Theorem with $c=0$, as $f(c)=0$ and $f'(c)=1$.
        To do so, we must calculate $f_k$ for the first few terms. Expanding $\sin_2^{-1}x$, we find
        \begin{equation*}
            f_0=0,\ f_1=1,\ f_2=0,\ f_3=1,\ f_4=0,\ f_5=9,\ f_6=0.
        \end{equation*}
        Using these values, we can find $\hat{f_k}$:
        \begin{equation*}
           \hat{f_1}=0, \  \hat{f}_2=\frac{1}{3},\ \hat{f}_3=0,\ \hat{f}_4=\frac{9}{5},\ \hat{f}_5=0.
        \end{equation*}
        We may now use these values to find the first few $g_n$ using the formulas given above. To this end, we record a couple of special Bell polynomials  that will be used below: $B_{n, n}(x_1) = (x_1)^n$ and $B_{n, n-1} (x_1, x_2) = {n \choose 2} (x_1)^{n-2}x_2$. These are obtained by simplifying the general Bell polynomial given above.
        
        When $n=1$,  $g_1=\frac{1}{f_1} = \frac{1}{1}=1$. When $n=2$, we have $g_2 = (-1)^1\cdot2^{(1)}B_{1,1}(0)=0$. Similarly, when $n=3$, we have
        \begin{eqnarray*}
g_3 & = & \frac{1}{f_1^3} \left(  (-1)^1 3^{(1)} B_{2,1} (\hat{f_1}, \hat{f_2}) + (-1)^2 3^{(2)} B_{2, 2}(\hat{f_1})\right) \\
& = &  \frac{1}{1^3} \left(  -3 B_{2,1} (0, 1/3) + 12 B_{2, 2}(0)\right)  \\
& = & -3 {2 \choose 2} \frac{1}{3} + 12 (0^2) = -1.
        \end{eqnarray*}

        In the same manner, applying this formula to the next few values of $n$, we find that $g_4=0$ and $g_5=1$.
    
        Substituting these values into the formula for $g(z)$ given by Lagrange Inversion Theorem above, we have:
        \begin{align*}
            g(z)&=0+\sum_{n=1}^\infty g_n \frac{(z-0)^n}{n!}. \\
            \sin_2(z) &=z-\frac{z^3}{3!}+\frac{z^5}{5!}+\cdots.
        \end{align*}
        When $p=4$, these computations get more tedious. Using SageMath, we find that
    \begin{equation*}
        \sin_4x=x-\frac{18}{5!}x^5+\frac{14364}{9!}x^9- \cdots.
    \end{equation*}
    \end{example}

The above ideas prove the following theorem.

\begin{thm} \label{thm:analytic}
 For any integer $p > 1$, the functions $\sin_p^{-1} x$ and $\sin_p x$ are analytic at $x = 0$.
\end{thm}

 It is well-known that $\sin x/x \rightarrow 1$ as $x \rightarrow 0$. We now generalize this result. 
 
 \begin{cor}
 Let $p >1$ be an integer. Then we have
 \[ \lim_{x \rightarrow 0} \frac{\sin_p x}{x} = 1.\]
 \end{cor}
 
 \begin{proof}
By Theorem \ref{thm:analytic}, we know that $\sin_p x$ is analytic at $x = 0$, and moreover, from the CIVP, $\sin_p 0 = 0$. Therefore, we can express $\sin_p x$ as a power series whose constant term is $0$:
 \[\sin_p x = a_1x + a_2x^2 + a_3 x^3 + \cdots + a_nx^n+ \cdots. \]
 Differentiating both sides and invoking the CIVP gives:
 \[(\cos_p x)^{p-1} = a_1+2a_2x + 3a_3x^2 + \cdots + nx^{n-1} + \cdots. \]
 Since $\cos_p(0) =1$, setting $x=0$ in the above equation  tells us that $a_1 = 1$.
 Finally, we have 
 \[\lim_{x \rightarrow 0} \, \frac{\sin_p x}{x} =  \lim_{x \rightarrow 0} \,  \frac{x+a_2x^2 + a_3x^3 + \cdots}{x} = \lim_{x \rightarrow 0} \, 1+a_2x +a_3x^2  + \cdots = 1.\]
 \end{proof}

Note that from this result it also follows that  $\tan_p x/x \rightarrow 1$ as $x \rightarrow 0$. One can also prove these limits using l'H\^{o}pital's  rule. In the same vein, one can also show the following.

 \begin{cor}
 \[ \lim_{x \rightarrow 0} \frac{\sin_4 x  - x}{x^5} = -\frac{18}{5!}.\]
 \end{cor}

\subsection{Rigidity} Note that the missing terms of the Taylor series for $\sin_4(x)$ are exactly the ones that were also missing in $\sin_4^{-1}(x)$; see Example \ref{talorexample}. In fact, for both functions, the non-zero terms in the Taylor series correspond to powers of $x$ that form an arithmetic progression of the form $4m+1$.  We proved this fact in Theorem \ref{thm1} for $\sin_n^{-1}(x)$. We now conjecture that this is also true for $\sin_n(x)$. 

\begin{conjecture}  \label{conj}
Let  $n$ be a positive integer.  Then 
\[ \left(\frac{d^l}{dx^l} \sin_n (x) \right) \biggr \rvert_{x=0}  \ne 0  \iff l \equiv 1 \mod n. \]
\end{conjecture}

This led to the following, more general question in analysis.
\vskip 3mm

\noindent
\textbf{Question:}  Suppose $f(x)$ is a real-valued function that is infinitely differentiable at $x=a$ such that  $f’(a) \ne 0$. Let $f(a) = b$ and let $g(x)$ be the local inverse of $f(x)$ (this exists  by the inverse function theorem) at $x = a$.   Is it true that for every positive integer $n$, the $n$th derivative of $f(x)$ at $x=a$ is non-zero if and only if the nth derivative of $g(x)$ at $x=b$ is non-zero?

\vskip 3mm
It turns out that, in general, the above answer is no. Take for example $f(x)=x^2$. We have $f(1) = 1$ and $f'(1) = 2 \ne 0$. At $x=1$, the local inverse of $f(x)$ is $g(x) = \sqrt{x}$.  Note that for all $k \ge 3$, $f^{(k)}(1) = 0$ but $g^{(k)}(1) \ne 0$. On the other hand, for the function  $f(x) = \sin(x)$, the above question has an affirmative answer because the Taylor series for $\sin x$ and $\sin^{-1} x$, have only odd terms.  This leads naturally to the following definition. 

\vskip 3mm
\begin{defn}
Let  $y=f(x)$ be a function that is infinitely differentiable at $x=a$ such that  $f’(a) \ne 0$.  We say that $f(x)$ is rigid at $x=a$ if  for any positive integer $k$, $f^{(k)}(a) \ne 0$ if and only if $g^{(k)} (b) \ne 0$, where $g(x)$ is the local inverse of $f(x)$ at $x = a$ and $b = f(a)$.
\end{defn}

In this terminology, $f(x) =x^2$ is not rigid at $x=1$ but $f(x) = \sin x$ is rigid at $x=0$. Conjecture \ref{conj} can now be restated as follows. For any positive integer $k$, $\sin_k x$ is rigid at $x=0$.

\vskip 3mm
\noindent
\textbf{Question:} What are necessary and sufficient conditions for a function $y= f(x)$ that is infinitely differentiable at $x=a$ to be rigid at $a$?

\section{Generalized  $\pi$ values} \label{pi_p}

\subsection{Organic definition} 
As we generalize trigonometric functions in the $p$-norm, we must also take into consideration generalizing the value of $\pi$. Recall that  $\pi =2\sin^{-1}(1)$. Using this as our inspiration, we can organically define $\pi_p$ as $\pi_p:=2\sin_p^{-1}(1)$. Using our $\sin_p^{-1}x$ formula we derived in Section \ref{Inverse functions} and letting $x=1$, we get 
\begin{equation}\label{pi_peq}
    \pi_p= 2\int_{0}^{1}{ \frac{1}{(1-t^p)^{\frac{p-1}{p}}}} \,dt.
    \end{equation}
Note that, unless otherwise indicated, when we refer to $\pi$, we are referring to  $\pi_2$.

When $p=2$, we find that the area of the unit circle is equal to $\pi$. It is then natural to wonder if $\pi_p$ has any relation to the area of a unit $p$-circle.

\begin{prop}
The area of a unit $p$-circle is $\pi_p$, when $p \ge 1$.
\end{prop}
\begin{proof}
In Proposition \ref{areapara}, we found that the area of the sector of the $p$-circle that connects the points $(x,y)$ and $(1,0)$ is given as a function of $y$ by $a(y)=\frac{1}{2}\int_{0}^{y}(1-t^p)^{\frac{1}{p}-1} \,dt$. If we let $(x,y)$ be the point $(0,1)$, we get the area of the unit $p$-circle in the first quadrant, given by $a(1)=\frac{1}{2}\int_{0}^{1}(1-t^p)^{\frac{1}{p}-1} \,dt$. Since the unit $p$-circle has 4-fold symmetry, we can multiply both sides of the equation by four to find the area of the entire $p$-circle:
\begin{equation*}
    4a(1)=2\int_{0}^{1}(1-t^p)^{\frac{1}{p}-1} \,dt.
\end{equation*}
From Section \ref{Inverse functions}, we know that $2\int_{0}^{x}(1-t^p)^{\frac{1}{p}-1} \,dt=2\sin_p^{-1}x$. Therefore, we know that the right hand side of the equation is $2\sin_p^{-1}(1)$, which is equal to $\pi_p$ as shown in Equation (\ref{pi_peq}). We have also already established that $a(1)$ is the area of the quarter unit $p$-circle, so $4\cdot a(1)$ gives us the area of the entire unit $p$-circle. Therefore, we find that the area of the unit $p$-circle is $\pi_p$. 
\end{proof}

\begin{cor}\label{0>pi_p>4}
For any $p \ge 1$, we have  $2 \le \pi_p <4$.
\end{cor}
\begin{proof}
As shown above, $\pi_p$ is the area of a unit $p$-circle. When $p = 1$, we get the region bounded by the square $|x| + |y| = 1$, which has area $2$.
 Similarly, since the $p$-circle is inscribed in a square of side length two, we know that the area of the $p$-circle is bounded by the square's area, which is 4. This shows that for any $p \ge 1$, we have  $2 \le \pi_p<4$.
\end{proof} 

\subsection{A formula for $\pi_p$}
We now show how we can compute $\pi_p$ in terms of the gamma function. To this end, we need another special function called  
the beta function, $\beta(x,y)$, which is closely related to the gamma function and can be defined as $\beta(x,y)=\int_{0}^{1} t^{x-1}(1-t)^{y-1} \,dt$, for any two real  numbers $x,y$ such that $x>0$ and $y>0$. We can put the beta function in terms of gamma using the property
\begin{equation}\label{BetaGamma}
    \beta(x,y)=\frac{\Gamma(x)\Gamma(y)}{\Gamma(x+y)}.
\end{equation}

\begin{prop} For any $p \ge 1$, we have
 \begin{equation*}
    \pi_p =\frac{2\Gamma^2(\frac{1}{p})}{p\Gamma(\frac{2}{p})}.
\end{equation*}
In particular, $\pi_p$ is a differentiable function of $p$.
\end{prop}

\begin{proof}
Referring back to our definition for $\pi_p$, if we let $u=t^p$, we can express $\pi_p$ in terms of the beta function  as follows:
\begin{equation*}
    \pi_p= 2\int_{0}^{1}{ \frac{1}{(1-t^p)^{\frac{p-1}{p}}}} \,dt = \frac{2}{p} \int_{0}^{1} (1-u)^{\frac{1}{p}-1}\cdot u^{\frac{1}{p}-1}\,du = 2 \beta(1/p, 1/p).
\end{equation*}
We can then use Equation \eqref{BetaGamma} to put $\pi_p$ in terms of the gamma function, and that gives the formula stated in the proposition.
Since $\Gamma(x)$ is a differentiable function and compositions and quotients of differentiable functions are again differentiable, it follows that $\pi_p$ is differentiable.
\end{proof}
\begin{example} \label{pi_p:2,3,4}
 Using the above equation, we can numerically approximate  $\pi_p$ for any $p$. For $p=2, 3$ and $4$, we get:

    \begin{equation*}
    \pi_2=\frac{2\Gamma^2(\frac{1}{2})}{2\Gamma(1)}\approx 3.1415,
       \quad 
       \pi_3=\frac{2\Gamma^2(\frac{1}{3})}{3\Gamma(\frac{2}{3})}\approx3.533
       \quad \text{and}\quad
       \pi_4=\frac{2\Gamma^2(\frac{1}{4})}{4\Gamma(\frac{2}{4})}\approx3.708.
    \end{equation*}
    \end{example}

\subsection{Properties of $\pi_p$}

We have already seen that $\pi_p$ is a differentiable function of $p$ for all $p > 0$. Is it a monotonic function? 
Example \ref{pi_p:2,3,4} suggests that $\pi_p$ increases with $p$. We now prove that fact.

\begin{prop}
  $\pi_p$ is an increasing function on $(0, \infty)$.
\end{prop}
\begin{proof}
Recall that $\pi_p$ is the area of a unit $p$-circle. Since $p$-circles have a 4-fold symmetry, we get  $\pi_p= 4\int_{0}^{1}{(1-x^p)^{\frac{1}{p}}} \,dt$. We will be done if we can show that, for any fixed value of $x$ in $(0, 1)$, $(1-x^p)^{\frac{1}{p}}$ is an increasing function in $p$.  This is because if 
$(1-x^{p_1})^{\frac{1}{p_1}} <  (1-x^{p_2})^{\frac{1}{p_2}}$ for all $x$ in $(0, 1)$ and $p_1 < p_2$, then 
\[ 4 \int_0^1 (1-x^{p_1})^{\frac{1}{p_1}} \, dx < 4 \int_0^1 (1-x^{p_2})^{\frac{1}{p_2}} \, dx,\]
showing that $\pi_{p_1} < \pi_{p_2}$ whenever 
$0 < p_1 < p_2$.

To this end, let $\psi(p):=(1-x^p)^{\frac{1}{p}}$, where $x$ is a fixed number in $(0,1)$. Taking the natural logarithm  and differentiating with  respect to $p$ on both sides,   we get
    \begin{eqnarray*}
     \ln(\psi(p)) &=& \frac{\ln(1-x^p)}{p}, \\
      \frac{\psi'(p)}{\psi(p)} & =&  \frac{-\ln(1-x^p)}{p^2}+\frac{-\ln(x)x^p}{p(1-x^p)}, \\
      \psi'(p) & =&  (1-x^p)^{\frac{1}{p}} \left( \frac{-\ln(1-x^p)}{p^2}+\frac{-\ln(x)x^p}{p(1-x^p)} \right).
    \end{eqnarray*}
For $0 < x < 1$ and $p >0$,  note that $0 <  1-x^p < 1$. Therefore, $\ln(x)$ and $\ln(1-x^p)$ are both negative. This shows that all parts of the derivative are positive. Therefore,  $\psi'(p) > 0$, which means $\psi(p)$ is an increasing function.
\end{proof}
Having shown above that $\pi_p$ is an increasing function and that it has an upper bound of 4 in Corollary \ref{0>pi_p>4}, we know that a limit exists. It is then only natural to wonder what the limit of $\pi_p$ is.
\begin{prop}
 $\lim_{p\to\infty} \pi_p=4$.
\end{prop}

\begin{proof}
 Using $\pi_p=\frac{2\Gamma^2(\frac{1}{p})}{p\Gamma(\frac{2}{p})}$, we can take the limit of $\pi_p$ as $p$ approaches infinity. Note that we have not yet stated Legendre's duplication formula, 
 ${\Gamma(2z)}=\frac{\Gamma(z)\Gamma(z+\frac{1}{2})}{2^{1-2z}\sqrt{\pi}}$.

\begin{align*}
    \pi_p&=\frac{2}{p}\cdot\frac{\Gamma^2(\frac{1}{p})}{\Gamma(\frac{2}{p})} \\
    &=\frac{2}{p}\cdot\frac{\Gamma^2(\frac{1}{p})\cdot2^{1-\frac{2}{p}}\sqrt{\pi}}{\Gamma(\frac{1}{p})\Gamma\left(\frac{1}{p}+\frac{1}{2}\right)} \\
    &=\frac{2}{p}\cdot\frac{\Gamma(\frac{1}{p})\cdot2^{1-\frac{2}{p}}\sqrt{\pi}}{\Gamma\left(\frac{1}{p}+\frac{1}{2}\right)} \\
     &=\frac{2\cdot \Gamma(\frac{1}{p}+1)\cdot2^{1-\frac{2}{p}}\sqrt{\pi}}{\Gamma\left(\frac{1}{p}+\frac{1}{2}\right)}.
       \end{align*}
  \[ \lim_{p\to\infty} \pi_p=\frac{2\cdot \Gamma(1)\cdot2\sqrt{\pi}}{\Gamma(\frac{1}{2})} 
    =\frac{ 1 \cdot 4\sqrt{\pi}}{\sqrt{\pi}} 
    =4. \]
\end{proof}

\begin{figure}[H]
    \centering
    \includegraphics[width=0.4\textwidth]{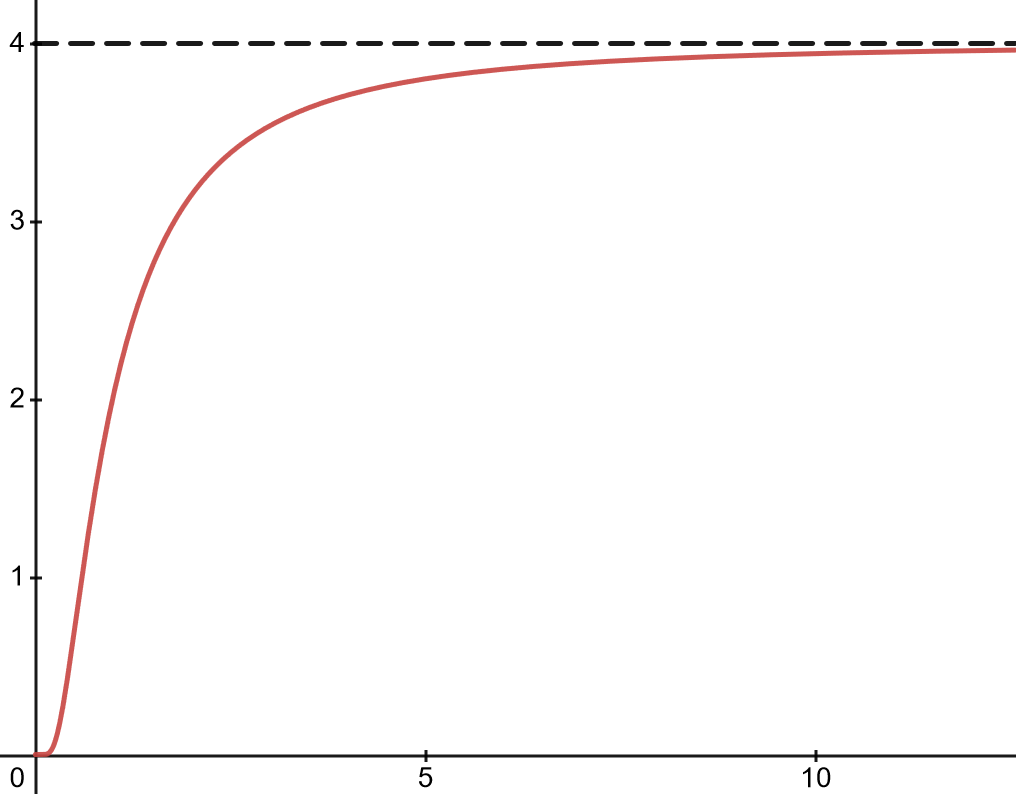}
    \caption{Graph of $\pi_p$ }
    \label{fig:pi_p}
\end{figure}

Recall that it is well-known that $\pi_2 = \pi$ is an irrational number (a number that is not the ratio of two integers). On the other hand, $\pi_1 = 2$, a rational number. ($\pi_1 =2$ because it is the area enclosed by the square $|x| +|y| = 1$ of side length $\sqrt{2}$.)  It is natural to ask for what values of $p$, is $\pi_p$ irrational? This is a hard question. Since $\pi_p$ is a continuous function, it takes rational and irrational values infinitely often; see Figure \ref{fig:pi_p}.

\subsection{$\pi_p$ with the Monte Carlo method}
    Since we have shown that $\pi_p$ can be described as the area of a $p$-circle, we can use a rather fun technique to approximate the value of $\pi_p$. Given a $p$-circle shaped dartboard inscribed inside a square, what is the probability that a uniformly random throw will land on the dartboard (assuming that the dart must land inside the square)? The probability is the ratio of the area of the board to the area of the box. Therefore, if we have $n$ throws where $t$ of them land on the dartboard, the probability would be $\frac{t}{n}=\frac{\pi_p}{4}$. We can solve for $\pi_p$ to get $\pi_p=\frac{4t}{n}$. Because of the law of large numbers, when $n\to\infty$, the ratio goes to the true ratio, and we find the true value of $\pi_p$. Writing a simple program to do this for us, at $n=1,000,000,000$, we get $\pi_3 \approx 3.53324$ and $\pi_4 \approx 3.7081$.

\section{Optimal unit $p$-circles} \label{optimal}
One question that naturally arises when examining $p$-circles is, ``At what value of $p$ is the corresponding squircle halfway between a unit circle ($p=2$) and a square ($p \to\infty$)?" This question was examined from three lenses: area, perimeter, and curvature. 
\subsection{Area}
We sought to find the value of $p$ for which the area enclosed by the $p$-circle $|x|^p+|y|^p=1$ is $\frac{\pi+4}{2}$, which is the average of the areas of the unit circle and the square that the unit circle is inscribed in. Because the $p$-circle is symmetric, we can examine the first quadrant only, resulting in the following equation:
\begin{equation*}
    \int_0^1\sqrt[p]{1-x^p}\ dx=\frac{\pi+4}{8}.
\end{equation*}
Using SageMath, the root of this equation can be found, giving the approximation $p\approx 3.162038$. As such, we can conclude that the value of $p$ for which the area of the $p$-circle is exactly halfway between the areas of the unit circle ($p=2$) and the square in which it is inscribed is $p\approx 3.162038$.

\subsection{Perimeter}
Next, we want to find the value of $p$ for which the perimeter of a unit $p$-circle is halfway between those of a unit circle and the square that the unit circle is inscribed in. The circumference of a unit circle is $2\pi$, and the perimeter of a square that contains the unit circle is 8. Therefore, we have to find the value of $p$ for which the perimeter of a unit $p$-circle is $\pi+4$. 
To find the perimeter of the unit $p$-cirle, we apply the Euclidean arc length formula to the defining equation of a $p$-circle. We equate the resulting integral to $\pi+4$ to obtain the equation:
\[
 \pi+4 = 4\int_{0}^{1} {\sqrt{1+(1-x^p)^{2(1-p)/p} x^{2(p-1)}} \; dx.}
 \]
 
We solved this equation numerically using SageMath to find that  $p \approx 4.667489$.

\subsection{Curvature} Finally, we wish to find $p$ such that the curvature of the unit $p$-circle is halfway between that of a square (here said to have curvature 0) and the $2$-unit circle (which has curvature $1$).
For a given smooth curve $C$ in $\mathbb{R}^2$, the curvature is
a measure of how different our curve is from a circle at a given point. While there are many equivalent
formulations of the curvature of a given curve, the following gives the curvature for a curve defined implicitly by 
$F(x,y) = 0$:
\begin{align*}
    \kappa = \frac{|F_y^2F_{xx} - 2F_xF_yF_{xy} + F_x^2F_{yy}|}{(F_x^2+F_y^2)^{\frac{3}{2}}}.
\end{align*}

Using the relation $F(x, y) = x^p + y^p - 1 = 0$ for the unit $p$-circle, we obtain
$$\kappa = (p-1)\frac{x^py^{2p}+x^{2p}y^p}{(x^{2p}y^2+y^{2p}x^2)^{\frac{3}{2}}}(xy).$$

If we investigate the curvature of the unit $p$-circle at the point $x = y$, we find that
\begin{align*}
    \kappa &= (p-1)x^2 \frac{2x^{3p}}{(2x^{2p+2})^{\frac{3}{2}}} = \frac{p-1}{\sqrt{2}x}.
\end{align*}

When $x = y$, we can write the relation for the unit circle as $2x^p = 1$, which gives $x = 2^{-\frac{1}{p}}$. Substituting for $x$ gives

\begin{align*}
    \kappa &=  \frac{p-1}{\sqrt{2} \cdot 2^{-\frac{1}{p}}} = \frac{p-1}{2^{\frac{1}{2}-\frac{1}{p}}} = (p-1)2^{\frac{1}{p}- \frac{1}{2}}.
\end{align*}

Therefore, if we solve for $p$ such that the unit $p$-circle has curvature 1/2, we find that $p \approx 1.43643264$.

\subsection{Resulting graphs}
\begin{figure}[H]
    \centering
    \includegraphics[width=\textwidth]{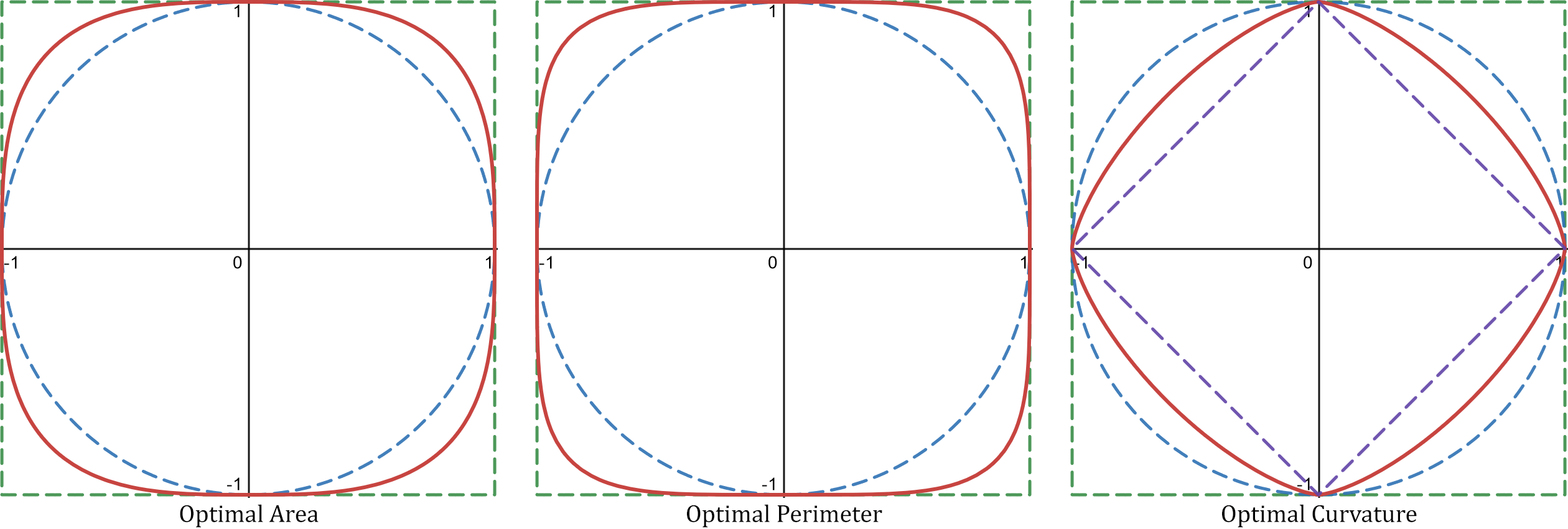}
    \caption{All 3 optimal $p$-circles}
    \label{fig:opgraphs}
\end{figure}

Graphing these 3 $p$-circles gives Figure \ref{fig:opgraphs} where the unit circle and square are dashed, and the $p$-circle is solid. For the optimal curvature, we also have $p=1$ since both $p=1$ and $p\to\infty$ have the same curvature.

\section{Rational Points on $p$-circles} \label{rational}
\subsection{2-circles and Pythagorean triples}
Right triangles (and as a result, Pythagorean triples) have long been objects of mathematical interest, studied intensely by the Babylonians even more than a thousand years before Pythagoras \cite{NO}. 
Given any Pythagorean triple $(x,y,z)$ satisfying $x^2 + y^2 = z^2$, we may divide all parts by $z^2$ to obtain 
$\frac{x^2}{z^2} + \frac{y^2}{z^2} = 1$. Then the point $(\frac{x}{z}, \frac{y}{z})$ is a rational point which lies on the $2$-unit circle defined by $x^2 + y^2 = 1$. On the other hand, given any rational number $\frac{v}{u}$, we obtain
the Pythagorean triple $(u^2-v^2, 2uv, u^2 + v^2)$ \cite{SJ}. In this manner, we may translate contexts between rational points on the unit circle and right triangles with integer side lengths.

\subsection{$p$-circles and Fermat's Last Theorem}
We can generalize the known results for rational points on $2$-circles and ask the same question for $p$-circles where $p$ is an integer greater than 2. There are certainly 4 trivial rational points along the axes of the graph, which are the points $(0,1), (1,0), (0,-1), (-1,0)$. To find the others, we may look at the rational solutions in the first quadrant  and use symmetry to extend our answers to the entire unit $p$-circle. 

Let $p$ be an integer greater than 2 and let $P = (\frac{p_1}{q_1},\frac{p_2}{q_2}) \in \mathbb{Q}^2$ be a rational point on the unit $p$-circle lying in the first quadrant. As $P$ lies on the unit $p$-circle and $P$ is in the first quadrant, we must have $(\frac{p_1}{q_1})^p + (\frac{p_2}{q_2})^p = 1$ and $\frac{p_1}{q_1}>0, \frac{p_2}{q_2}>0$. We then find 
$(p_1q_2)^p + (p_2q_1)^p = (q_1q_2)^p$. However, by Fermat's Last Theorem, there are no positive integers $p_1, p_2, q_1, q_2$ that satisfy this relation. Thus, there exist no rational solutions in the first quadrant. Then, by symmetry, we see that the only rational points on the circle are exactly those along the axes.

\section{Future Research} \label{conclusion}
The results in this paper seem to indicate that $p$-trigonometric functions have interesting but complex behavior. For instance, even basic formulas such as the double-angle formulas for $\sin (2x)$  and $\cos (2x)$ seem to have no straightforward generalization. Similarly, understanding higher derivatives of $\sin_p(x)$ at $x=0$ looks very difficult; see Conjecture \ref{conj} and the questions following it.   There are several other open questions. We list a few that we think merit further study.

\begin{enumerate}
    \item We know the derivatives of $\sin_p x$ and $\cos_p x$. What about $ \int \sin_p x \, dx $ and $\int \cos_p x \, dx$?  Using the Taylor series for $\sin_p x$ and $\cos_p x$, one can evaluate these integrals as series. But are there closed-form answers for these integrals? 

    \item The parametrization of $p$-circles we considered in this paper are with respect to area. We can also parametrize these curves with respect to arc length.  These give yet another generalization of the $p$-trigonometric functions.  What properties do these functions have? 
    
    \item Can we extend this work for $(p, q)$-trigonometric functions that come from looking at the curves $|x|^p + |y|^q = 1$? Parametrizing these curves will give us $\sin_{p, q} x$ and $\cos_{p, q} x$. What can be said about these functions?
   
    \item So far, we have been working in  $\mathbb{R}^2$. Can we extend this work to $\mathbb{R}^3$? To this end, we should look at the unit $p$-sphere $|x|^p + |y|^p +|z|^p= 1$. For $p=2$, this is the standard unit sphere, and as $p$ goes to infinity, we get a cube that encloses the unit sphere. These surfaces can be called sphubes ($p$-spheres), analogous to our squircles ($p$-circles). It opens gates to a whole new area of research. What are the parametric equations of these surfaces? Can we do sphubical trigonometry that is similar to spherical trigonometry? What are the volume and surface areas of the regions enclosed by these surfaces? What is the Gaussian curvature function of these surfaces? 
\end{enumerate}

\begingroup
\raggedright

\bibliographystyle{plain}
\bibliography{citations}
\endgroup

\end{document}